\documentclass{ws-procs9x6}
\input{liemacs-red.sty} 

\newcommand{\bV}{{\mathbb V}}
 
\newcommand\cO{\mathcal{O}}

\newcommand{\derat}[1]{\frac{d}{dt} \hbox{\vrule width0.5pt
                height 3mm depth 3mm${{}\atop{{}\atop{\scriptstyle t=#1}}}$}}

\newtheorem{nremark}[theorem]{Remark}
\newtheorem{nproposition}[theorem]{Proposition}

\newtheorem{ndefinition}[theorem]{Definition}
\newtheorem{nlemma}[theorem]{Lemma}
\newtheorem{nexample}[theorem]{Example}

\begin{document}

\title{Semi-bounded unitary representations of infinite-dimensional 
Lie groups} 

\author{K.-H. Neeb} 

\address{Fachbereich Mathematik, TU Darmstadt, \\
Schlossgartenstrasse 7, 64285 Darmstadt, Germany \\ 
$^*$E-mail: neeb@mathematik.tu-darmstadt.de}

\begin{abstract} In this note we introduce the concept of a semi-bounded 
unitary representations of an infinite-dimensional Lie group $G$. 
Semi-boundedness 
is defined in terms of the corresponding momentum set in the dual 
$\g'$ of the Lie algebra $\g$ of $G$. After dealing with some 
functional analytic issues concerning certain weak-$*$-locally compact 
subsets of dual spaces, called semi-equicontinuous, 
we characterize unitary representations which are bounded in the sense
that their momentum set is  equicontinuous, we characterize 
semi-bounded representations of locally convex spaces in terms 
of spectral measures, and we also describe a method to compute 
momentum sets of unitary representations of reproducing kernel 
Hilbert spaces of holomorphic functions.  
\end{abstract} 


\bodymatter

\section*{Introduction}  \label{sec:0}

For infinite-dimensional Lie groups it seems quite hopeless to develop 
a general theory of unitary representations. This is mainly due to two 
obstacles which are non-existent for finite-dimensional Lie groups. 
The first one is that there is no group algebra such as 
$L^1(G)$ or $C^*(G)$, whose representations 
are in one-to-one correspondence with the continuous unitary representations 
of $G$. The second one is that there is no general structure theory 
for infinite-dimensional Lie groups, such as the Levi decomposition and 
the fine structure theory of semisimple Lie groups. 

However, there are many interesting classes of infinite-dimensional 
Lie groups which possess a rich unitary representation theory. 
Many of these representations show up naturally in various contexts 
of mathematical physics (\refcite{Mick87, Mick89, PS86, Go04, SeG81, Se58, Se78}), 
where the Lie algebra $\g = \L(G)$ 
of the group under consideration often contains an element $h$, corresponding 
to the Hamiltonian of the underlying physical system, for which the 
spectrum of the operator $i\cdot \dd \pi(h)$ 
in the ``physically relevant'' representations $(\pi, \cH)$ 
is bounded from below. This suggests to study 
representations of infinite-dimensional Lie groups in 
terms of semi-boundedness properties of spectra. 

Let $G$ be a Lie group in the category of smooth manifolds 
modelled on locally convex spaces for which a smooth 
exponential function $\exp_G \: \g \to G$ exists
\begin{footnote}{We shall keep this assumption throughout this paper.}  
\end{footnote}
(cf.\ \refcite{Mil84, Ne06, GN08}).  
For a unitary representation $(\pi,\cH)$ of $G$ we write 
$\pi^v(g) := \pi(g)v$ for its orbits maps 
and call the representation $(\pi, \cH)$ {\it smooth} 
if the space 
$$\cH^\infty := \{ v \in \cH \: \pi^v \in C^\infty(G,\cH)\}$$ 
of smooth vectors is dense in $\cH$. 
Then all operators $i \cdot \dd \pi(x)$, $x \in \g$, are 
essentially selfadjoint \cite[Lemma~3.6]{Ne08} and crucial 
information on their spectrum is contained in the momentum set $I_\pi$ 
of the representation, a subset of 
the topological dual $\g'$ of the Lie algebra $\g$ of $G$. 
It is defined as the weak-$*$-closed convex hull of the image of the 
momentum map on the projective space of $\cH^\infty$
$$ \Phi_\pi \: \bP({\cal H}^\infty)\to \g' \quad \hbox{ with } \quad 
\Phi_\pi([v])(x) 
= \frac{1}{i}  \frac{\la  \dd\pi(x).v, v \ra}{\la v, v \ra}\quad \mbox{ for } 
[v] = \C v. $$
As a weak-$*$-closed convex subset, the momentum set is completely determined 
by its support functional 
\begin{equation}
  \label{eq:supp-spec}
s_\pi \: \g \to \R \cup \{\infty\}, \quad s_\pi(x) 
= - \inf \la I_\pi,x \ra = \sup(\Spec(i \cdot \dd \pi(x)))
\end{equation}
(cf.\ \cite[Lemma~3.7]{Ne08}). 

It is now natural to study those representations for which 
$s_\pi$, resp., the set $I_\pi$, contains the most significant 
information. A natural regularity condition is that 
the function $s_\pi$ is bounded on some non-empty open subset of $\g$. 
We call such representations {\it semi-bounded}. 
Then the domain $s_\pi^{-1}(\R)$ of $s_\pi$ is a convex cone 
with non-empty interior and $s_\pi$ is continuous on this 
open cone (cf.\ \cite[Prop.~4.8]{Ne08}). 
Since the momentum set $I_\pi$ is invariant under the coadjoint action, 
the function $s_\pi$ and its domain are invariant under the adjoint action. 
In \refcite{Ne08} we described a framework which permits us to find 
$C^*$-algebras whose representations are in one-to-one correspondence 
with certain classes of semi-bounded representations of~$G$.

In the present note we discuss several issues concerning semibounded 
representations. In Section~\ref{sec:1} we briefly discuss the 
functional analytic background for semi-bounded representations. 
In particular we introduce the concept of a semi-equicontinuous set 
of linear functionals in such a way that semi-bounded representations 
are those whose momentum set is semi-equicontinuous. 
In many concrete situations, it is desirable to calculate 
momentum sets of representations by restricting the momentum map 
to quite small invariant submanifolds of projective space. 
The main result of Section~\ref{sec:2} is a theorem 
which describes criteria under which the momentum set 
of a representation on a reproducing kernel 
Hilbert space of holomorphic functions 
can be determined directly from the reproducing kernel. 
Section~\ref{sec:3} is devoted to a characterization of bounded 
representations as those for which the corresponding homomorphism 
$\pi \: G \to U(\cH)$ is a smooth homomorphism of Lie groups, where the unitary 
group $U(\cH)$ carries its natural Banach--Lie group structure. 
This result is a quite straight forward consequence of the definitions 
if $G$ is locally exponential, but for a general Lie group we have 
to use some more sophisticated arguments. 
Finally we explain in Section~\ref{sec:4} that for the case where 
$G = (E,+)$ is the additive group of a locally convex space, 
semi-bounded representations are precisely those obtained by 
functional calculus from spectral measures supported by semi-equicontinuous 
subsets of the dual space. 

For the future, the most important 
issue concerning semi-bounded representations 
is to understand the structural implications of the existence of 
semi-bounded representations, to determine which invariant 
semi-equi\-continuous subset of $\g'$ actually arise as momentum sets 
of unitary representations and to understand their convexity properties 
and the related symplectic geometry. 

\section{Semi-equicontinuous convex sets} \label{sec:1}

Let $E$ be a real locally convex space and $E'$ its topological dual, i.e., 
the space of continuous linear functionals on~$E$. 
We write $\la \alpha, v \ra = \alpha(v)$ for the natural pairing 
$E' \times E \to \R$ and endow $E'$ with the weak-$*$-topology, i.e., 
the coarsest topology for which all linear maps 
$$ \eta_v \: E' \to \R, \quad \eta_v(\alpha) := \alpha(v) $$
are continuous. For a subset $X \subeq E'$, the set 
$$ B(X) := \{ v \in E \: \inf\la X, v \ra > - \infty\} $$
is a convex cone which coincides with the domain 
of the {\it support function} 
$$ s_X \: E \to \R \cup \{ \infty\}, \quad s_X(v) 
:= - \inf \la X, v \ra = \sup \la X, - v\ra $$
of $X$ in the sense that $B(X) = s_X^{-1}(\R)$. As a sup 
of a family of continuous linear functionals, the function 
$s_X$ is convex, lower semicontinuous and positively homogeneous. 

\begin{nremark}\label{rem:1.1} (a) The 
set $X$ is weak-$*$-bounded if and only if all functions 
$\eta_v$ are bounded on $X$, i.e., $B(X) = E$. It is equicontinuous 
if and only if the function $s_X$ is bounded on some $0$-neighborhood in 
$X$. Each equicontinuous subset is in particular weak-$*$-bounded, 
but the converse only holds if $E$ is a barrelled space 
(Uniform Boundedness Principle) 
\cite[Ch.~III, \S 4, no.~2, th.]{Bou07}. 
Recall that all Fr\'echet spaces are barrelled, but that also  
locally convex direct limits of barrelled spaces are barrelled.

(b) If $Y := \oline{\conv}(X)$ denotes the weak-$*$-closed convex 
hull of $X$, then $B(X) = B(Y)$, $s_X = s_Y$, 
and, using the Hahn--Banach Separation Theorem, 
$Y$ can be reconstructed from $s_Y$ 
by 
$$Y = \{ \alpha \in E' \: (\forall v \in B(Y))\, 
\alpha(v) \geq \inf \la Y,v \ra = - s_Y(v)\}.$$
\end{nremark}

\begin{ndefinition} We call $X$ {\it semi-equicontinuous} if 
$s_X$ is bounded on some non-empty open subset of $E$. 
\end{ndefinition} 

If $X$ is semi-equicontinuous, then $B(X)$ clearly has interior points. 
We also have a partial converse: 
\begin{nproposition} \label{prop:1.3}\cite[Thm.~4.10]{Ne08} 
If $E$ is barrelled, then a subset $X \subeq E'$ is 
semi-equicontinuous if and only if $B(X)$ has interior points. 
\end{nproposition}

\begin{nproposition} \label{prop:1.4}\cite[Props.~4.13, 4.4]{Ne08} If $X$ is a semi-equicontinuous 
weak-$*$-closed convex subset, then the following assertions hold: 
\begin{description}
\item[\rm(a)] For each $v\in B(X)^0$, 
$\eta_v \: X \to \R$ is a proper function. 
In particular, $X$ is locally compact in the weak-$*$-topology. 
\item[\rm(b)] $X = \{ \alpha \in E' \: (\forall v \in B(X)^0)\, 
\alpha(v) \geq  \inf \la X, v \ra \}.$
\end{description}
\end{nproposition}

\begin{nexample} \label{ex:1.5} (a) If $W \subeq E$ is an open convex cone and 
$$ W^\star := \{ \alpha \in E' \: \alpha(W) \subeq \R_+ \} $$
its dual cone, then $B(W^\star) = \oline W$ and $s_{W^\star} = 0$ on $W$. 
Therefore $W^\star$ is semi-equicontinuous. 

(b) Let $Y$ be a topological space and $\omega \: Y \to ]0,\infty[$ a 
non-zero continuous function. Then 
$$ C_\omega(Y,\R) := \Big\{ 
f \in C(Y,\R) \: \sup\frac{|f|}{\omega} < \infty \Big\} $$
is a Banach space with respect to the norm $\|f\|:= \sup\frac{|f|}{\omega}$. 
Each element $y \in Y$ defines a continuous linear functional 
on this space by $\delta_y(f) := f(y)$, and the set 
$$ X := \{ \delta_y \: y \in Y \} \subeq C_\omega(Y,\R)' $$ 
is semi-equicontinuous. In fact, $\omega$ is positive on $Y$, 
and the open unit ball $B_1(\omega)$ around $\omega$ consists of non-negative 
functions. Therefore $X$ is contained in the dual of an open cone, 
hence semi-equicontinuous. 
\end{nexample}

\begin{nremark} (a) If $\phi \: E \to F$ is a continuous linear map 
between locally convex spaces and $X \subeq F'$ is semi-equicontinuous, 
then the adjoint map $\phi' \: F' \to E', \alpha \mapsto \alpha \circ \phi$, 
maps $X$ into a semi-equicontinuous subset because 
$s_{\phi'(X)} = s_X \circ \phi$ is bounded on some non-empty open 
subset of $E$. 

(b) If $X \subeq E'$ is semi-equicontinuous and 
$v \in B(X)^0$, then $\eta_v$ is bounded from below, so that 
for some $c > 0$, the function $\omega := \eta_v + c$ is positive 
on $X$. For any other $w \in E$ there exists an $\eps > 0$ 
with $v \pm \eps w \subeq B(X)$, so that there exists a $d > c$ with 
$\eta_{v\pm \eps w}(X) \geq -d$. This implies that 
$|\eta_w| \leq \eps^{-1}(\eta_v + d)$ on $X$, which in turn implies that 
$\eta_w \in C_\omega(X,\R)$. We thus obtain a map 
$$ \eta \: E \to C_\omega(X,\R), \quad w \mapsto \eta_w\res_{X} $$
which is easily seen to be continuous with $\eta'(\delta_\alpha) = \alpha$ 
for each $\alpha \in X$. 

This observation shows that any semi-equicontinuous 
set is contained in the image 
of $\{ \delta_x \:x \in X\}$ under the adjoint of some continuous linear map 
$\phi \: E \to C_\omega(X,\R)$, and, in view of (a) and 
Example~\ref{ex:1.5}, we know that, conversely, all such 
sets are semi-equicontinuous. 
\end{nremark}

\section{Momentum sets of smooth unitary representations} \label{sec:2}

Let $G$ be a Lie group with a smooth exponential 
function $\exp_G \: \g \to G$, i.e., all the curves 
$\gamma_x(t) := \exp_G(tx)$, $x \in \g$, are one-parameter groups with 
$\gamma_x'(0) = x$. 

Now let $(\pi, \cH)$ be a smooth unitary representation 
of $G$. We  then obtain for each $x \in\g$ 
a unitary one-parameter group $\pi_x(t) := \pi(\exp_G(tx))$ and 
for each $v \in \cH^\infty$ the derivative  
$\dd \pi(x)v := \derat0 \pi(\exp_G(tx))v$
exists and defines an unbounded operator $\dd \pi(x) \: \cH^\infty 
\to \cH^\infty \subeq \cH$. Since the map $\cH^\infty \to C^\infty(G,\cH), 
v \mapsto \pi^v$ is equivariant with respect to the 
representation of $G$ on $C^\infty(G,\cH)$ by $(g.f)(x) := f(xg)$, 
it easily follows that 
$$ \dd \pi \: \g \to \End(\cH^\infty) $$
defines a representation of $\g$, called the {\it 
derived representation}. 

We also have a natural action of $G$ on the projective 
space $\bP(\cH^\infty)$ by $g.[v] := [\pi(g)v]$ and the coadjoint 
action of $G$ on $\g'$ by $\Ad^*(g).\alpha := \alpha \circ \Ad(g)^{-1}$. 
An easy calculation now shows that the momentum map 
$$ \Phi_\pi \: \bP(\cH^\infty) \to \g', \quad 
\Phi_\pi([v])(x) = \frac{\la \dd\pi(x)v,v\ra}{i\la v,v\ra} $$
is $G$-equivariant. In particular, its image is $G$-invariant, and therefore 
$I_\pi$ is invariant under the coadjoint action. 

\begin{nproposition} \label{prop:2.1} \cite[Lemma~3.7]{Ne08} For each $x \in \g$, the 
closure of the operator $\dd \pi(x)$ on $\cH^\infty$ is the 
infinitesimal generator of the unitary one-parameter group 
$\pi_x(t) := \pi(\exp_G(tx))$ and satisfies 
$$ \sup(\Spec(i\cdot \dd\pi(x))) = s_{I_\pi}(x) = - \inf \la I_\pi,x\ra. $$
\end{nproposition}

In terms of the momentum set, defined in the introduction, 
we now define: 

\begin{ndefinition}\label{def:2.2} A smooth representation 
$(\pi,\cH)$ is called {\it bounded}, resp., {\it semi-bounded} if 
its momentum set $I_\pi \subeq \g'$ is 
equicontinuous, resp., semi-equicontinuous. 
\end{ndefinition}

\begin{nremark} (a) In view of Proposition~\ref{prop:2.1}, the convex cone 
$B(I_\pi)$ is the set of all elements of $x$ for which the 
selfadjoint operator $i \dd\pi(x)$ is bounded above, which in  
turn is equivalent to the existence of an extension to a semigroup 
homomorphism
$$ \hat\pi_x \: \C_+ := \R + i \R_+ \to B(\cH), \quad 
\hat\pi_x(z) := e^{z \dd\pi(x)} $$
(where the exponential is to be understood in terms of the functional 
calculus with respect to a spectral measure) which is strongly continuous and holomorphic on the open upper half plane 
(cf. \refcite{Ne00}). 

(b) If $(\pi, \cH)$ is a bounded representation for which 
$\ker(\dd\pi) = I_\pi^\bot = \{0\}$, then $\|\dd\pi(x)\|$ 
defines a $G$-invariant norm on $\g$. If $G$ is finite-dimensional, 
then the existence of an invariant norm implies that the 
Lie algebra $\g$ is compact. 

(c) If $\g$ is infinite-dimensional, then the existence of an invariant 
norm does not imply that the Lie bracket extends to the corresponding Banach 
completion. A simple example is the Lie algebra 
$(C^\infty(\T^2,\R),\{\cdot,\cdot\})$ of smooth functions 
on the $2$-torus, endowed with the Poisson bracket with respect to the 
canonical symplectic form $\omega = \dd x \wedge \dd y$. Then the 
$L^2$-inner product 
$$ (f,g) := \int_{\T^2} fg \cdot \omega $$ 
is invariant under the adjoint 
action of the corresponding Lie group $\Ham(\T^2, \omega)$ (which is 
simply given by translation), but the Poisson--Lie bracket 
$$ \{f,g\} = 
\frac{\partial f}{\partial x} \frac{\partial g}{\partial y} 
-\frac{\partial f}{\partial y} \frac{\partial g}{\partial x} $$
is not continuous with respect to the $L^2$-inner product. 
\end{nremark}

\begin{nremark} Semi-bounded unitary representations of finite-dimensional 
Lie groups have been studied in some detail in \refcite{Ne00}, where 
it is shown that all these representations are direct integrals 
of irreducible semi-bounded representations 
\cite[Sect.~XI.6]{Ne00}.  
and that, 
on the Lie algebra level, the irreducible representations 
are highest weight representations $(\pi_\lambda,\cH_\lambda)$  
\cite[Thms.~X.3.9, XI.4.5]{Ne00}.   
If $[v_\lambda] \in \bP(\cH^\infty)$ is a highest weight vector, then 
the corresponding $G$-orbit $G.[v_\lambda]$ has the remarkable 
property that 
$$ \Phi(G.[v_\lambda]) = \Ext(I_{\pi_\lambda}) = \Ad^*(G)\Phi([v_\lambda]) 
\quad \mbox{ and } \quad 
I_{\pi_\lambda} = \conv(\Phi(G.[v_\lambda])) $$ 
\cite[Thm.~X.4.1]{Ne00}. 
Moreover, two irreducible semi-bounded representations are equivalent if 
and only the corresponding momentum sets, resp., the coadjoint orbits 
of extreme points coincide \cite[Thm.~X.4.2]{Ne00}.  
\end{nremark}

One major feature of unitary highest weight representations 
is that the image of the highest weight orbit already determines the 
momentum set as the closed convex hull of its image. 
It is therefore desirable to understand in which situations 
smaller subsets of $\bP(\cH^\infty)$ already determine the momentum set. 
As we shall see below, this situation frequently occurs when 
$\cH$ consists of holomorphic functions on some complex manifold. 

\begin{ndefinition} Let $M$ be a complex manifold (modelled on a locally 
convex space) and $\cO(M)$ the space of holomorphic complex-valued 
functions on $M$. We write $\oline M$ for the conjugate complex manifold. 
A holomorphic function 
$$ K \: M \times \oline M \to \C $$
is said to be a {\it reproducing kernel} of a Hilbert subspace 
$\cH \subeq \cO(M)$ if for each $w\in M$ the function 
$K_w(z) := K(z,w)$ is contained in $\cH$ and satisfies 
$$ \la f, K_z \ra = f(z) \quad \mbox{ for } \quad z \in M, f \in \cH. $$
Then $\cH$ is called the a {\it reproducing kernel Hilbert space} 
and since it is determined uniquely by the kernel $K$, it is 
denoted $\cH_K$ (cf.\ \cite[Sect.~I.1]{Ne00}). 

Now let $G$ be a real Lie group and 
$\sigma \: M \times G \to M, (m,g) \mapsto m.g$ 
be a smooth right action of $G$ on $M$ by holomorphic maps. Then 
 $(g.f)(m) := f(m.g)$ defines a unitary representation 
of $G$ on a reproducing kernel Hilbert space $\cH_K \subeq \cO(M)$ 
if and only if the kernel $K$ is invariant: 
$$ K(z.g, w.g) = K(z,w) \quad \mbox{ for } \quad z,w \in M, g \in G. $$
In this case we call 
${\cal H}_K$ a $G$-invariant reproducing kernel Hilbert space 
and write $(\pi_K(g)f)(z) := f(z.g)$ for the corresponding 
unitary representation of $G$ on $\cH_K$. 
\end{ndefinition} 

\begin{nlemma} Let $G$ be a Fr\'echet--Lie group, 
$(\pi_K,\cH_K)$ a unitary representation on a $G$-invariant reproducing kernel 
Hilbert space in $\cO(M)$ and 
$\Omega := \{ [K_m] \: m \in M, K(m,m)> 0\}$. 
Then the following assertions hold: 
\begin{description}
\item[(a)] $K_z \in \cH^\infty$ for each $z \in M$ and 
the representation of $G$ on $\cH_K$ is smooth. 
\item[(b)] If $x \in \g$ is such that 
the smooth action $M \times \R \to M, \break (m,t) \mapsto m.\exp_G(tx)$ 
extends to a holomorphic action of the upper half plane, then  
$$ \inf \la I_\pi, x\ra = \inf \la \Phi_{\pi_K}(\Omega), x\ra. $$
\end{description}
\end{nlemma} 

\begin{proof} (a) For each $f \in \cH$ and $z \in M$ we have 
$$ \la f, g.K_z \ra = \la g^{-1}.f, K_z \ra = (g^{-1}.f)(z) 
= f(z.g^{-1}), $$
which is a smooth function $G \to \C$. Hence the map 
$\alpha_z \: G \to \cH, g \mapsto g.K_z$, is weakly smooth. 
This implies that for each smooth map $h \: \R^n \to G$ the 
composition $\alpha_z \circ h$ is weakly smooth, 
hence smooth by Grothendieck's Theorem (cf.\ \cite[p.484]{Wa72}). 
Now we apply \cite[Rem.~12.5]{BGN04} 
to see that $\alpha_z$ is smooth. This means that 
each $K_z$ is a smooth vector, and since these elements 
span a dense subspace of $\cH$, (a) follows. 

(b) Let $(m,s) \mapsto m.s$ denote the holomorphic action of 
$\C_+$ on $M$ extending the given action of $\R$. 
For $s \in \C_+$ we put $s^* := - \oline s$, which turns $\C_+$ 
into an involutive semigroup. 
For $z,w \in M$, the functions 
$$f_1(s) := K(z.s,w) \quad \mbox{ and } \quad f_2(s) := K(z,w.s^*).$$
Both are holomorphic on $\C_+^0$, continuous on $\C_+$ and 
coincide on $\R$, so that they are equal \cite[Lemma~A.3.6]{Ne00}. 
On the dense subspace 
$$\cH_K^0 := \Spann \{ K_z \: z \in M \}$$ of 
$\cH_K$ we now obtain a representation of $\C_+$ by 
$(\hat\pi_x(s).f)(m) := f(m.s)$
(cf.\ \cite[Prop.~II.4.3]{Ne00}). 
Next we observe that 
\begin{align*}
& \frac{1}{2} \derat0 K(m.it,m.it)
= \frac{1}{2} \derat0 K(m.2it,m)\\
&= i \derat0 K(m.t,m)
= i \derat0 K_m(m.t)= i \derat0 (\exp_G(tx).K_m)(m)\\
&= i \derat0 \la \pi_K(\exp_G(tx)).K_m, K_m \ra 
= i\la \dd\pi_K(x)K_m, K_m \ra. 
\end{align*}
Therefore \cite[Prop.~IV.1]{Ne00b} implies that 
\begin{equation}
  \label{eq:2.2}
\|\hat\pi_x(a+ib)\| = e^{b\sup \la \Phi_{\pi_K}(\Omega),-x\ra}. 
\end{equation}
Clearly, 
$\sup \la I_{\pi_K}, -x \ra \geq 
 \sup \la \Phi_{\pi_K}(\Omega), -x \ra$ and 
if the right hand side is infinite, then both are equal. Suppose that 
this is not the case, so that $\hat\pi_x$ actually defines a 
representation of $\C_+$ by bounded operators on $\cH_K$. 

For each $f \in \cH_K$ and $z,w \in M$ we then have 
$\hat\pi_x(s).K_m = K_{m.s^*}$ 
\cite[Prop.~II.4.3]{Ne00}, so that 
$\la \hat\pi_x(s).K_m, f\ra 
=  \la K_{m.s^*}, f \ra = \oline{f(m.s^*)}$
is holomorphic on $\C_+^0$ and continuous on $\C_+$. 
Since the representation $\hat\pi_x$ on $\C_+$ is locally bounded, 
\cite[Lemma~IV.2.2]{Ne00} implies that 
$\hat\pi_x \: \C_+ \to B(\cH_K)$ is strongly continuous and 
holomorphic on $\C_+^0$. Now \cite[Lemma~VI.5.2]{Ne00} 
shows that 
$\hat\pi_x(s) = e^{s \dd \pi(x)}$ for each $s \in \C_+$, 
so that $i\dd\pi(x)$ is bounded above with 
$$ \|\hat\pi_x(i)\| = e^{\sup\Spec(i\dd\pi(x))} 
= e^{\sup \la I_{\pi_K}, -x \ra}.$$
Comparing with \eqref{eq:2.2} now completes the proof. 
\end{proof}

\begin{theorem}
  \label{thm:2.7} 
Let $G$ be a Fr\'echet--Lie group acting smoothly by holomorphic 
maps on the complex manifold $M$ and 
$\cH_K \subeq \cO(M)$ be a $G$-invariant reproducing kernel 
Hilbert space. If, for each $x \in B(I_{\pi_K})^0$,  
the action $(m,t) \mapsto m.\exp_G(tx)$ 
of $\R$ on $M$ extends to a holomorphic action of the upper half plane 
$\C_+$, then 
$$ I_{\pi_K} = \oline{\conv}(\Phi(\{ [K_m] \: K(m,m) = \|K_m\|^2 > 0\})). $$
\end{theorem}

\begin{proof} From the previous lemma, we obtain for each 
$x \in B(I_{\pi_K}^0)$ the relation 
$\inf \la I_\pi, x\ra = \inf \la \Phi_{\pi_K}(\Omega), x \ra$, 
so that the theorem follows from the reconstruction 
formula Proposition~\ref{prop:1.4}(b). 
\end{proof}

\begin{nremark} Since a holomorphic section of a vector bundle 
$\bV \to M$ can always be identified with a holomorphic function 
on the total space of the dual bundle $\bV' \to M$, any 
reproducing kernel Hilbert space of holomorphic sections 
can be realized as a reproducing kernel Hilbert space of holomorphic 
functions. Therefore the preceding theorem also applies to this 
more general situation. 
\end{nremark}

\section{Bounded representations} \label{sec:3}

The main goal of this section is to prove the following theorem 
characterizing bounded representations. The main difficulty of 
the proof is to bridge the gap between the smoothness of a unitary 
representation as defined above and the smoothness of an action 
of $G$ on the whole Hilbert space $\cH$. 

\begin{theorem}\label{thm:contin} If $(\pi,\cH)$ is a smooth representation 
of the Lie group $G$ with exponential function, 
then $I_\pi$ is equicontinuous, i.e., $\pi$ is bounded, 
if and only if 
$\pi \: G \to U(\cH)$ is a morphism of Lie groups, where 
$U(\cH)$ carries its natural Banach--Lie group structure. 
\end{theorem}

\begin{nremark} \label{rem:3.1} 
For $x \in \g$, the condition that 
$I_\pi(x)$ is a bounded subset of $\R$ means that 
the unitary one-parameter group 
$\pi_x(t) := \pi(\exp_G tx) = e^{t \dd \pi(x)}$
is norm-continuous (cf.\ \cite[Lemma~3.7]{Ne08}). 
This is the special case $G = \R$ of Theorem~\ref{thm:contin}. 
\end{nremark}

In the following we write $x.g \in T_g(G)$ for the $g$-right 
translate of $x \in \g \cong T_\1(G)$. 

\begin{nlemma} \label{lem:unique} 
Let $\sigma \: G \times M \to M, (g,m) \mapsto g.m$, be a smooth action 
of the connected Lie group $G$ on the smooth manifold $M$ and 
$$\dot\sigma \: \g \to \cV(M), \quad 
\dot\sigma(x)(m) := T_{(\1,m)}(\sigma)(-x,0) $$
the corresponding derived homomorphism of Lie algebras. 
If $f \: G \to M$ is a smooth map with 
$$ f(\1) = m \quad \mbox{ and } \quad 
T_g(f)(x.g) = -\dot\sigma(x)(f(g))\quad \mbox{ for } \quad g \in G, x\in \g, $$
then $f(g) = g.m$ for each $g \in G$, i.e., 
$f$ is the orbit map of $m$. 
\end{nlemma} 

\begin{proof} We consider the smooth map 
$h \: G \to M, h(g) := g^{-1}.f(g) = \sigma_{g^{-1}}(f(g))$
and calculate for $x\in \g$: 
\begin{align*}
T_g(h)(xg) 
&= T_{f(g)}(\sigma_g^{-1})T_g(f)(x.g)  + T_{f(g)}(\sigma_g^{-1})
\dot\sigma(x)(f(g)) \\
&= -T_{f(g)}(\sigma_g^{-1})\dot\sigma(x)(f(g)) + T_{f(g)}(\sigma_g^{-1})
\dot\sigma(x)(f(g)) =0.
\end{align*}
Since $G$ is connected, this 
implies that $h$ is constant $h(\1) = m$, which implies 
the lemma. 
\end{proof}

\begin{nproposition} \label{prop:unique-rep} Let 
$\pi_i \: G \to \GL(E)$, $i =1,2$, be two representations of the 
connected Lie group $G$ on the locally convex space $E$. 
We assume that $\pi_1$ is smooth and that for $\pi_2$ all orbit maps 
of element in the dense subspace $E^\infty \subeq E$ are smooth, so 
that we obtain a homomorphism of Lie algebras
$$ \dd \pi_2 \: \g \to \End(E^\infty), \quad 
\dd\pi_2(x)(v) := T_\1(\pi_2^v)(x) \quad \mbox{ for }
 \quad \pi_2^v(g) := \pi_2(g)v. $$
If these two representations are compatible in the sense that 
$$ \L(\pi_1)(x)v = \dd\pi_2(x)v \quad \mbox{ for } \quad v \in E^\infty, $$
then $\pi_1 = \pi_2$. 
\end{nproposition}

\begin{proof} We consider $\pi_1$ as a smooth action 
$\sigma(g)(v) := \pi_1(g)v$ of $G$ on $E$. Then the corresponding 
homomorphism 
$\dot\sigma \: \g \to {\cal V}(V)$ is given by 
$$ \dot\sigma(x)(v) = - \L(\pi_1)(x)(v). $$

Next, let $v \in E^\infty$. For the smooth map 
$f \: G \to E, f(g) := \pi_2(g)v$, we then have 
$$ T_g(f)(x.g) = \dd\pi_2(x)\pi_2(g)v = \L(\pi_1)(x)(\pi_2(g)v), $$
so that Lemma~\ref{lem:unique} implies that 
$\pi_2(g)v = f(g) = \pi_1(g)v$ for each $v \in E^\infty$. 
Since $\pi_1(g)$ and $\pi_2(g)$ are continuous operators on $E$ and 
$E^\infty$ is dense, it follows that $\pi_1(g) = \pi_2(g)$ for 
each $g \in G$. 
\end{proof}

\begin{proof}[Proof of Theorem~\ref{thm:contin}] Suppose first that $\pi$ is a morphism of Lie groups. 
Then ${\cal H}^\infty = {\cal H}$ and 
$\L(\pi) = \dd\pi\: \g \to \fu(\cH) = \L(U(\cH))$ is a continuous 
linear map. Hence $x \mapsto \|\dd\pi(x)\|$ defines a continuous 
seminorm on $\g$. For each $v \in {\cal H}$ we now have 
$|\la \dd\pi(x)v,v\ra| \leq \|\dd\pi(x)\| \cdot \|v\|^2,$
hence $|\la I_\pi, x \ra| \leq \|\dd\pi(x)\|$. Therefore 
$I_\pi$ is equicontinuous. 

Suppose, conversely, that $I_\pi$ is equicontinuous. 
Then  \cite[Lemma~3.7]{Ne08} implies that the essentially skew-adjoint 
operator $\dd\pi(x)$ on ${\cal H}^\infty$ extends to a bounded 
operator, also denoted by $\dd\pi(x)$, on $\cH$. 
Since the map \break 
$\dd\pi \: \g \to \End(\cH^\infty)$ is a representation of 
$\g$, the density of $\cH^\infty$ in $\cH$ and the estimate 
$$ \|\dd\pi(x)v\| \leq (\sup|\la I_\pi, x \ra|) \cdot \|v\| $$
imply that the map $\dd \pi \: \g \to B(\cH)$ also is linear, 
continuous and a homomorphism of Lie algebras. 

Since homomorphisms of Lie groups are smooth if and  only if 
they are 
smooth in an identity neighborhood, we may w.l.o.g.\ 
assume that $G$ is connected. 
Let $q_G \: \tilde G \to G$ be the simply connected covering 
group of $G$. 
Since the Banach--Lie group $U(\cH)$ is regular, the morphism 
$\dd\pi$ of Lie algebras integrates to a smooth group homomorphism 
$$ \tilde\pi \: \tilde G \to U(\cH)\quad \mbox{ with }\quad 
\L(\tilde\pi) = \dd\pi $$ 
(\refcite{Mil84}).
From Proposition~\ref{prop:unique-rep} we now derive that 
$\tilde\pi = \pi\circ q_G$, and this implies that 
$\pi$ is smooth. 
\end{proof}

If a unitary representation $\pi \: G \to U(\cH)$ is a morphism of Lie 
groups, then it is in particular norm continuous, i.e., a morphism 
of topological groups. One may now ask under which circumstances, the 
norm continuity implies that $\pi$ is smooth. 

\begin{nproposition}\label{prop:contin} Let $(\pi,\cH)$ be a smooth 
unitary representation of the Lie group $G$ 
which is norm continuous. Then $\pi$ is a morphism of 
Lie groups if $G$ is locally exponential or 
$\g$ is barrelled. 
\end{nproposition}

\begin{proof} (a) If $G$ is locally exponential, 
then $G$ and $U(\cH)$ are locally exponential Lie groups, 
and the smoothness of any continuous homomorphism follows from the 
Automatic Smoothness Theorem \cite[Thm.~IV.1.18]{Ne06}. 

(b) Suppose that $\g$ is barrelled. 
For each $x\in \g$, the unitary representation 
$\pi_x(t) := \pi(\exp_G(tx))$ is norm continuous, hence 
$I_\pi(x)$ is bounded. This implies that $I_\pi \subeq \g'$ is 
weak-$*$-bounded, and since $\g$ is barrelled, it is equicontinuous 
(Remark~\ref{rem:1.1}). Now Theorem~\ref{thm:contin} shows that 
$\pi$ is a morphism of Lie groups. 
\end{proof}

\section{The abelian case} \label{sec:4}

Let $G := (E,+)$ be a locally convex space, considered as a Lie group. 
We fix a weak-$*$-closed convex semi-equicontinuous 
subset $X \subeq E' = \g'$ and recall 
from Proposition~\ref{prop:1.4} that $X$ is locally compact w.r.t.\ 
the weak-$*$-topology.  
The following theorem characterizes the semibounded 
smooth unitary representations of $G$ with $I_\pi \subeq X$: 

\begin{theorem} For a smooth representation $(\pi,\cH)$ of 
$(E,+)$, the following are equivalent: 
\begin{description}
\item[\rm(a)] $I_\pi \subeq X$.  
\item[\rm(b)] There exists a holomorphic non-degenerate representation 
$$ \hat\pi \: S := E + i B(X)^0 \to B(\cH) $$
of involutive semigroups (with respect to $(x+iy)^* := -x + iy$), 
satisfying $\pi(v)\hat\pi(s) = \hat\pi(v + s)$ 
for $v \in E$, $s \in S$ and 
$$ \|\hat\pi(x+iy)\| \leq e^{- \inf \la X, y \ra}. $$
\item[\rm(c)] There exists a Borel spectral measure $P$ on the locally 
compact space $X$ with $P(X) = \1$ and $P(e^{i\eta_v}) = \pi(v)$ for 
each $v \in E$. 
\end{description}
\end{theorem}

\begin{proof} (a) $\Rarrow$ (b): For each $z = x + iy \in S$ we define 
$$ \hat\pi(x+iy) := \pi(x) e^{i\dd\pi(y)} \in B(\cH). $$
Since the spectral measures of the unitary one-parameter groups 
of the form $t \mapsto \pi(tv)$ commute pairwise, it follows that 
$\hat\pi$ is a homomorphism of semigroups compatible with the 
involution. From 
$$ \|\hat\pi(x+iy)\| = e^{\sup\Spec(i\dd\pi(y))} 
= e^{- \inf \la I_\pi, y\ra} $$
it follows that $\hat\pi$ is locally bounded. 
For each finite-dimensional subspace $F \subeq E$ intersecting 
$B(X)^0$ non-trivially, we apply \cite[Section~VI.5]{Ne00} 
to see that $\hat\pi$ is holomorphic on $F + i (F \cap B(X)^0)$, 
hence holomorphic since it is locally bounded (\refcite{He89}). 

(b) $\Rarrow$ (c) In view of 
\cite[Thm.~5.1]{Ne08}, the representation 
of $S$ comes from a representation of the $C^*$-algebra 
$C_0(X)$ with respect to the homomorphism 
$\gamma \: S \to C_0(X),\gamma(s)(\alpha) := e^{i\alpha(s)}$. 
Since the representations of $C_0(X)$ are in one-to-one 
correspondence with spectral measures on $X$ and $\pi$ is 
uniquely determined by its compatibility with $\hat\pi$, 
(c) follows. 

(c) $\Rarrow$ (a): Since $\eta_v$ is bounded from below on 
$X$ by $- s_X(v)$, we derive from (c) that 
$\sup(\Spec(i\dd\pi(v))) \leq s_X(v)$, which implies (a) 
(Proposition~\ref{prop:1.4}(b)). 
\end{proof}

\end{document}